\documentclass[10pt,reqno]{amsart}

\usepackage[utf8]{inputenc}
\usepackage[T1]{fontenc}
\usepackage{amsmath,amssymb,amsthm}
\usepackage{enumitem}
\usepackage{tikz}
\usepackage{xcolor}
\usepackage[margin=1.15in]{geometry}
\usepackage[hidelinks]{hyperref}

\newtheorem{theorem}{Theorem}[section]
\newtheorem{lemma}[theorem]{Lemma}
\newtheorem{proposition}[theorem]{Proposition}
\newtheorem{corollary}[theorem]{Corollary}
\theoremstyle{definition}
\newtheorem{definition}[theorem]{Definition}
\newtheorem{example}[theorem]{Example}
\theoremstyle{remark}
\newtheorem{remark}[theorem]{Remark}

\newcommand{\Z}{\mathbb{Z}}
\newcommand{\Fn}{F_2(C_n)}
\newcommand{\Kn}{C_n\boxtimes C_n}
\newcommand{\cint}[2]{[#1,#2]_n}
\newcommand{\lv}[2]{\langle #1,#2\rangle}
\newcommand{\Sset}[1]{\mathcal K_{#1}}
\newcommand{\Sn}{\Sset{n}}
\DeclareMathOperator{\dist}{d}

\begin{document}

\title[Independent sets in $F_2(C_n)$ and kings on the torus]{Independent sets in the
2-token graph of a cycle and kings on the torus}

\author{Luis Manuel Rivera}
\address{Unidad Acad\'emica de Matem\'aticas, Universidad Aut\'onoma de Zacatecas, Zacatecas, M\'exico}
\email{luismanuel.rivera@gmail.com}


\begin{abstract}
The $2$-token graph of the cycle $C_n$ and the king graph on the $n\times n$ torus are known to have the same independence number. Answering a question of Reyes, Dalf\'o and Fiol, for every $n\ge4$ we give an explicit bijection between a maximum independent set of the king graph and a maximum independent set of the $2$-token graph: each pair of tokens is represented by a king placed in the row of one of its two tokens and in the column given by their sum. We also show that the $2$-token graph of $C_n$ has exactly one maximum independent set, or exactly $n$ when $n\equiv3\pmod 4$, all of them rotations of each other. The argument is self-contained and yields a new
proof of the value of the independence number. 
\end{abstract}
\maketitle
\section{Introduction}\label{sec:intro}

Let $G$ be a graph and let $1\le m\le |V(G)|-1$. The \emph{$m$-token graph} $F_m(G)$ is the graph whose vertices are the $m$-subsets of $V(G)$, two of them being adjacent whenever their symmetric difference is an edge of $G$~\cite{FFHH}. For $m=2$, these graphs were studied earlier as \emph{double vertex graphs}~\cite{alavi}, and they are also known as \emph{symmetric squares}~\cite{alzaga,aude}. The \emph{strong product} of graphs $G$ and $H$, denoted by $G\boxtimes H$, is the graph with vertex set $V(G)\times V(H)$, in which two distinct vertices $(g,h)$ and $(g',h')$ are adjacent if and only if $g$ and $g'$ are equal or adjacent in $G$, and $h$ and $h'$ are equal or adjacent in $H$~\cite{HIK}. For $n\ge3$, the graph $\Kn$ is the \emph{king graph} of the $n\times n$ toroidal chessboard: its vertices are the squares of the board, and two squares are adjacent if and only if a king can move from one to the other in a single move. In particular, the independent sets of $\Kn$ are precisely the sets of squares occupied by pairwise non-attacking kings. By results of Hales~\cite{hales} and of Sonnemann and Krafft~\cite{Sonnemann} (see also~\cite{vesel,watkins}), and by a theorem of de Alba, Carballosa, Lea\~nos and Rivera~\cite{alba}, for every $n\ge4$,
\begin{equation}\label{eq:alphas}
\alpha(\Kn)=\left\lfloor\frac n2\left\lfloor\frac n2\right\rfloor\right\rfloor=\alpha(\Fn).
\end{equation}
For $n=3$ the equality also holds trivially, since $C_3\boxtimes C_3\cong K_9$ and $F_2(C_3)\cong C_3$.

The coincidence~\eqref{eq:alphas} was already noticed in~\cite{alba}, where it is observed that the values of $\alpha(\Fn)$ form the sequence A189889 of the OEIS~\cite{oeis}, which counts the maximum number of non-attacking kings on an $n\times n$ toroidal board. More recently, Reyes, Dalf\'o and Fiol~\cite{RDF} recalled this coincidence and remarked that it would be interesting to find a correspondence between the independent sets of $\Kn$ and those of $\Fn$. Theorem~\ref{thm:introB} and Corollary~\ref{cor:all} below provide such a correspondence. For the independence numbers of other $2$-token graphs and of pair graphs, which are the supertoken graphs $\mathcal{F}_2(G)$ in the terminology of~\cite{RDF}, see~\cite{JR}, and for $2$-token graphs of join graphs, see~\cite{RV}.
\subsection{Main results}

Throughout the rest of the paper, unless otherwise stated, $n\ge3$ and $k=\lfloor n/2\rfloor$. For $1\le\delta\le k$, let $D_\delta$ denote the set of vertices $\{u,v\}$ of $\Fn$ such that $\dist(u,v)=\delta$, where $\dist(\cdot,\cdot)$ denotes the distance in $C_n$. We call $D_\delta$ the \emph{$\delta$-th level} of $\Fn$. Since the distance between two distinct vertices of $C_n$ lies in $\{1,\dots,k\}$, the levels $D_1,\dots,D_k$ partition $V(\Fn)$. 

Our first result describes all the maximum independent sets of $\Fn$. Note that all the levels are independent sets, with one exception: when $n$ is odd, the top level $D_k$ induces a cycle of length $n$ (Proposition~\ref{prop:levels}), and then only part of it can be used. Below, a union of levels indexed by an empty range is understood to be empty, which occurs only for $k=1$.

\begin{theorem}\label{thm:introA}
Let $n\ge3$ and let $I$ be a set of vertices of $\Fn$. Then $I$ is a maximum independent set of $\Fn$ if and only if
\[
I=\begin{cases}
D_1\cup D_3\cup\cdots\cup D_{k-1}, & k\text{ even},\\[2pt]
D_1\cup D_3\cup\cdots\cup D_{k-2}\cup W, & k\text{ odd},
\end{cases}
\]
where $W=D_k$ if $n$ is even, and $W$ is a maximum independent set of the cycle induced by $D_k$ if $n$ is odd. In particular, $\alpha(\Fn)=\lfloor n\lfloor n/2\rfloor/2\rfloor$.
\end{theorem}
Let $I_n$ be the maximum independent set of $\Fn$ constructed in~\cite{alba}, which we recall in Definition~\ref{def:In}. Theorem~\ref{thm:introA} shows that every maximum independent set of $\Fn$ is a rotation of $I_n$:
\begin{corollary}\label{cor:introA}
For every $n\ge3$, the maximum independent sets of $\Fn$ are precisely the rotations $I_n+r$ with $r\in\Z_n$. There is only one of them when $n\not\equiv3\pmod4$, and there are exactly $n$ when $n\equiv3\pmod4$.
\end{corollary}
Both results are proved in Section~\ref{sec:classification}. The argument relies solely on the structure of $\Fn$ as a union of the levels $D_1,\dots, D_k$; in particular, it does not use the value of $\alpha(\Fn)$ obtained in~\cite{alba}, and hence it yields a self-contained proof of that value. For paths, the answer is simpler: $\alpha(F_2(P_n))=\lfloor n^2/4\rfloor$ and $F_2(P_n)$ has a unique maximum independent set for every $n\ge2$, namely the union of its odd levels, that is, the set of pairs of vertices at odd distance in $P_n$; see Remark~\ref{rem:paths}.

Our second result is a bijective account of the coincidence~\eqref{eq:alphas}. Identifying $V(C_n)$ with $\Z_n$ (see Section~\ref{sec:prelim}), let
\[
\ell_0=\{(x,2x):x\in\Z_n\}
\]
and consider the \emph{sum map}
\[
\psi\colon(\Z_n\times\Z_n)\setminus\ell_0\longrightarrow V(\Fn),\qquad \psi(x, y)=\{x, y-x\},
\]
which represents a pair of tokens $\{u,v\}$ by a king in column $u+v$ and in the row of one of the two tokens (more details are given in Section~\ref{sec:main}). The more obvious map $\pi(x,y)=\{x,y\}$, suggested by the isomorphism $\Fn\cong\bigl((C_n\square C_n)-\Delta\bigr)/S_2$, works for even $n$ but not for odd $n$; the sum map differs from it by the change of coordinates $(x,y)\mapsto(x,y-x)$ (see Section~\ref{sec:main}).

\begin{theorem}\label{thm:introB}
Let $n\ge4$. There is a maximum independent set $\Sn$ of $\Kn$ such that
$\Sn\cap\ell_0=\emptyset$ and $\psi|_{\Sn}\colon\Sn\to I_n$ is a bijection.
\end{theorem}

The set $\Sn$ is given explicitly in Section~\ref{sec:main} (Theorem~\ref{thm:main}),
as a union of pieces of lines of slope~$2$. Equivalently, one can choose, for each pair $\{u,v\}\in I_n$, one of its two tokens $r$ in such a way that the kings placed at the squares $(r,u+v)$, one for each pair, are pairwise non-attacking. Such a choice is given by the map $\Phi_n$ of Definition~\ref{def:Phi}, described by an explicit rule in Proposition~\ref{prop:Phi}; it is the inverse of the bijection $\psi|_{\Sn}$. A key property of the sum map is that it matches natural independent sets of $\Kn$ with the pieces of the maximum independent sets of $\Fn$: for $1\le\delta<n/2$, each of the two lines of slope $2$ and offset $\pm\delta$ in $\Z_n\times\Z_n$ is an independent set of $\Kn$ of size $n$, and $\psi$ maps it bijectively onto the level $D_\delta$ (Proposition~\ref{prop:lines}). By Theorem~\ref{thm:introA}, the levels $D_\delta$ with $\delta$ odd are precisely the pieces from which the maximum independent sets of $\Fn$ are built. Since $\psi$ is equivariant with respect to rotations, Corollary~\ref{cor:introA} and Theorem~\ref{thm:introB} together imply that every maximum independent set of $\Fn$ is obtained in this manner from a maximum independent set of $\Kn$ (Corollary~\ref{cor:all}).

We point out that Theorem~\ref{thm:introB} does not reprove~\eqref{eq:alphas}. Since $\Sn$ is an independent set of $\Kn$ with $|\Sn|=|I_n|=\alpha(\Fn)$, the theorem gives a constructive proof of the inequality $\alpha(\Kn)\ge\alpha(\Fn)$. Combined with Theorem~\ref{thm:introA}, this yields the lower bound $\alpha(\Kn)\ge\lfloor n\lfloor n/2\rfloor/2\rfloor$, so that the only ingredient of~\eqref{eq:alphas} that we take from~\cite{hales,Sonnemann} is the reverse inequality. What Theorem~\ref{thm:introB} adds is an explicit correspondence between the extremal configurations of the two graphs.

The coincidence~\eqref{eq:alphas} itself turns out to be exclusive to one and two tokens: in Section~\ref{sec:three} we show that $\alpha(F_m(C_n))<\alpha(C_n^{\boxtimes m})$ for every $n\ge4$ and every $3\le m\le n-1$, where $C_n^{\boxtimes m}$ is the king graph of the $m$-dimensional torus of side $n$.

\section{Preliminaries}\label{sec:prelim}
The vertex set of the cycle $C_n$ is identified with $\Z_n$, whose elements are represented by $1,2,\dots,n$; in particular, $0$ and $n$ denote the same element of $\Z_n$, and we use whichever representative is more convenient. Two vertices $u,v$ are adjacent if and only if $u-v\in\{1,-1\}$. Let
\[
\alpha_n:=\left\lfloor\frac n2\left\lfloor\frac n2\right\rfloor\right\rfloor=
\begin{cases}
k^2, & n=2k,\\[2pt]
nk/2, & n=2k+1,\ k\text{ even},\\[2pt]
(nk-1)/2, & n=2k+1,\ k\text{ odd}.
\end{cases}
\]
For $A\subseteq\Z_n$ and $r\in\Z_n$, we define $A+r=\{a+r:a\in A\}$. Note that $(\Z_n\setminus A)+r=\Z_n\setminus(A+r)$: since $x\mapsto x+r$ is a bijection of $\Z_n$, we have $y\in(\Z_n\setminus A)+r\iff y-r\notin A\iff y\notin A+r$. For a set $\mathcal A$ of $2$-subsets of $\Z_n$, we write $\mathcal A+r=\{A+r:A\in\mathcal A\}$. The assignment $A\mapsto A+r$ induces an automorphism of $\Fn$, which we refer to as a \emph{rotation}. For integers $a\le b$ with $b-a\le n-1$, the \emph{cyclic interval} $\cint ab\subseteq\Z_n$ is the set of residues of $a,a+1,\dots,b$; it contains $b-a+1$ elements and depends only on the residue of $a$ and on the difference $b-a$. Clearly, if $a\le b<c$ and $c-a\le n-1$, then $\cint ab$ and $\cint{b+1}c$ are disjoint; if, moreover, $c-a=n-1$, they partition $\Z_n$.

In the strong product $\Kn$, two distinct vertices $(x,y)$ and $(x',y')$ are adjacent (or \emph{attack} each other) if and only if $x'-x\in\{0,1,-1\}$ and $y'-y\in\{0,1,-1\}$. We refer to $x$ as the \emph{row} and to $y$ as the \emph{column} of the vertex $(x,y)$.

For a vertex $A=\{u,v\}$ of $\Fn$, let $\dist(A)=\dist(u,v)\in\{1,\dots,k\}$ denote the distance between $u$ and $v$ in $C_n$. For $c\in\Z_n$ we also write $|c|=\dist(0,c)$; thus $0\le|c|\le n/2$, and $|c|$ is the ordinary absolute value of any integer representative of $c$ lying in $[-n/2,n/2]$. For $x\in\Z_n$ and $1\le\delta\le k$, we write
\[
\lv x\delta=\{x,\,x+\delta\},
\]
so that the level $D_\delta$ defined in Section~\ref{sec:intro} is $D_\delta=\{A\in V(\Fn):\dist(A)=\delta\}=\{\lv x\delta:x\in\Z_n\}$. If $\delta<n/2$, the relation $\lv x\delta=\lv{x'}\delta$ implies $x=x'$, so $|D_\delta|=n$; if $n=2k$, then $\lv xk=\lv{x+k}k$ and $|D_k|=k$. 

We now recall the sets introduced in~\cite{alba}. For $1\le i\le n-1$, let
$L_i=\{\{j,\,n-i+j\}:1\le j\le i\}$; thus $L_i$ consists of the pairs $\{u,v\}$
with $1\le u<v\le n$ and $v-u=n-i$; every such pair has distance $\min\{i,n-i\}$;
and the sets $L_1,\dots,L_{n-1}$ partition $V(\Fn)$. Hence $L_\delta\cup L_{n-\delta}=D_\delta$
for $1\le\delta<n/2$ and $L_k=D_k$ when $n=2k$, while for $n=2k+1$ the identity
$\{j,j+k+1\}=\lv{j+k+1}k$ gives
\[
L_k=\{\lv xk:x\in\cint{k+2}{2k+1}\}\subseteq D_k,\qquad |L_k|=k.
\]
The maximum independent set of $\Fn$ constructed in~\cite{alba} is the union of those $L_i$ whose index $i$ is one of the integers $\delta$, $n-\delta$ with $\delta$ odd and $\delta\le k-1$, together with $L_k$ when $k$ is odd; in terms of the levels, it
can be described as follows.

\begin{definition}\label{def:In}
We define the set
\[
I_n=\begin{cases}
D_1\cup D_3\cup\cdots\cup D_{k-1}, & k\text{ even},\\[2pt]
D_1\cup D_3\cup\cdots\cup D_{k-2}\cup D_k, & n=2k,\ k\text{ odd},\\[2pt]
D_1\cup D_3\cup\cdots\cup D_{k-2}\cup L_k, & n=2k+1,\ k\text{ odd}.
\end{cases}
\]
\end{definition}
\begin{remark}\label{rem:In}
In the first two cases, $I_n$ is the union of all the levels of odd index. If $n$ is even, then the parity of an element of $\Z_n$ is well-defined and $\dist(u,v)\equiv v-u\equiv u+v\pmod 2$, so $I_n$ is precisely the set of pairs $\{u,v\}$ with $u+v$ odd, that is, with one odd and one even element.
\end{remark}

\section{\texorpdfstring{Proof of Theorem~\ref{thm:introA}}{Proof of Theorem 1.1}}\label{sec:classification}

We start by describing the adjacency of $\Fn$ in the coordinates $\lv x\delta$.

\begin{proposition}\label{prop:levels}
Let $x\in\Z_n$ and $1\le\delta\le k$. The neighbours of $\lv x\delta$ in $\Fn$ are
precisely the vertices listed in those of the following items that apply:
\begin{enumerate}[label=\textup{(\roman*)}]
\item $\lv{x}{\delta-1}$ and $\lv{x+1}{\delta-1}$, if $\delta\ge2$;
\item $\lv{x}{\delta+1}$ and $\lv{x-1}{\delta+1}$, if $\delta\le k-1$;
\item $\lv{x+k}{k}$ and $\lv{x-k}{k}$, if $n=2k+1$ and $\delta=k$;
\item $\lv{x+k}{k-1}$ and $\lv{x+k+1}{k-1}$, if $n=2k$ and $\delta=k$.
\end{enumerate}
Consequently, every edge of $\Fn$ joins two consecutive levels, except when $n=2k+1$, in which case $D_k$ induces the circulant graph on $\Z_n$ with connection set $\{k,-k\}$, that is, a cycle of length $n$. In particular, each level is independent, except for $D_k$ when $n$ is odd.
\end{proposition}

\begin{proof}
The neighbours of $A=\{x,x+\delta\}$ are the $2$-subsets among $\{x\pm1,x+\delta\}$ and $\{x,x+\delta\pm1\}$. If $\delta\ge2$, then $\{x,x+\delta-1\}=\lv x{\delta-1}$ and $\{x+1,x+\delta\}=\lv{x+1}{\delta-1}$, while for $\delta=1$ these two sets are not $2$-subsets. If $\delta\le k-1$, then $\{x,x+\delta+1\}=\lv x{\delta+1}$ and $\{x-1,x+\delta\}=\lv{x-1}{\delta+1}$. If $\delta=k$ and $n=2k+1$, then $\{x,x+k+1\}=\lv{x+k+1}{k}=\lv{x-k}k$ and $\{x-1,x+k\}=\lv{x+k}k$, because $(x+k+1)+k=x$ and $(x+k)+k=x-1$. If $\delta=k$ and $n=2k$, then $\{x,x+k+1\}=\lv{x+k+1}{k-1}$ and $\{x-1,x+k\}=\lv{x+k}{k-1}$.

For the last statement, note that in case (iii) the vertex $\lv xk$ is joined to $\lv{x\pm k}k$, and that the circulant graph on $\Z_n$ with connection set $\{\pm k\}$ is connected, hence a cycle of length $n$, because $\gcd(k,2k+1)=1$.
\end{proof}

Given an independent set $I$ of $\Fn$, we write $a_\delta=|I\cap D_\delta|$.

\begin{lemma}\label{lem:pair}
Let $1\le\delta\le k-1$ and assume $|D_{\delta+1}|=n$, that is, $n$ is odd or $\delta+1\le k-1$. Then the bipartite graph formed by the edges of $\Fn$ between $D_\delta$ and $D_{\delta+1}$ is a cycle of length $2n$. Consequently, for every independent set $I$ of $\Fn$,
\[
a_\delta+a_{\delta+1}\le n,
\]
and if equality holds, then $I$ contains one of the two levels and is disjoint from the other.
\end{lemma}

\begin{proof}
By Proposition~\ref{prop:levels}, the neighbours of $\lv x\delta$ in $D_{\delta+1}$ are $\lv x{\delta+1}$ and $\lv{x-1}{\delta+1}$, and the neighbours of $\lv y{\delta+1}$ in $D_\delta$ are $\lv y\delta$ and $\lv{y+1}\delta$. Hence the resulting bipartite graph $H$ is $2$-regular with $2n$ vertices, and it is connected because the walk
\[
\lv x\delta,\ \lv x{\delta+1},\ \lv{x+1}\delta,\ \lv{x+1}{\delta+1},\ \lv{x+2}\delta,\dots
\]
visits all of them. Thus $H$ is a cycle of length $2n$. Now $I\cap(D_\delta\cup D_{\delta+1})$ is an independent set of $H$, so $a_\delta+a_{\delta+1}\le\alpha(C_{2n})=n$. A cycle of even length has exactly two maximum independent sets, namely its two colour classes, which here are $D_\delta$ and $D_{\delta+1}$; so equality forces $I\cap(D_\delta\cup D_{\delta+1})$ to be $D_\delta$ or $D_{\delta+1}$.
\end{proof}

\begin{lemma}\label{lem:topeven}
Let $n=2k$ with $k\ge2$. Then $a_{k-1}+a_k\le n$ for every independent set $I$ of $\Fn$, and if equality holds, then $D_{k-1}\subseteq I$ and $I\cap D_k=\emptyset$.
\end{lemma}

\begin{proof}
For $x\in\Z_n$, let
\[
T_x:=\bigl\{\lv xk,\ \lv x{k-1},\ \lv{x+k}{k-1}\bigr\}.
\]
Since $\lv xk=\lv{x+k}k$, we have $T_x=T_{x+k}$, and the sets $T_1,\dots,T_k$ partition $D_{k-1}\cup D_k$: each of the $k$ vertices of $D_k$ lies in exactly one of them, and the $2k$ vertices of $D_{k-1}$ are distributed in pairs $\{\lv x{k-1},\lv{x+k}{k-1}\}$. Applying Proposition~\ref{prop:levels}(i) to both representations $\lv xk$ and $\lv{x+k}k$ of the same vertex, we see that it is adjacent to both $\lv x{k-1}$ and $\lv{x+k}{k-1}$, while $\lv x{k-1}$ and $\lv{x+k}{k-1}$ are non-adjacent because $D_{k-1}$ is independent. Hence $T_x$ induces a path on three vertices, so $|I\cap T_x|\le2$, with equality if and only if $\lv x{k-1},\lv{x+k}{k-1}\in I$ and $\lv xk\notin I$. Summing over $x=1,\dots,k$ gives $a_{k-1}+a_k\le 2k=n$, with equality if and only if $D_{k-1}\subseteq I$, in which case $I\cap D_k=\emptyset$.
\end{proof}

\begin{proof}[Proof of Theorem~\ref{thm:introA}]
We first check that each of the sets described in the statement is independent and has $\alpha_n$ elements. By Proposition~\ref{prop:levels}, all the edges of $\Fn$ join consecutive levels, except for the edges inside $D_k$ when $n$ is odd. In such a set, the levels that occur have pairwise non-consecutive indices; each of them is independent (Proposition~\ref{prop:levels}) and, when $n$ and $k$ are both odd, $W$ is independent by assumption; hence the set is independent. We now compute their cardinality. All of them have the same cardinality as $I_n$, since in the last case $|W|=\alpha(C_n)=k=|L_k|$. If $k$ is even, then $I_n$ is the union of the $k/2$ levels $D_1,D_3,\dots,D_{k-1}$, each of size $n$, so $|I_n|=\frac k2\,n=\alpha_n$. If $n=2k$ with $k$ odd, then $|I_n|=\frac{k-1}2\,n+|D_k|=k(k-1)+k=k^2=\alpha_n$. Finally, if $n=2k+1$ with $k$ odd, then $|I_n|=\frac{k-1}2\,n+k=\frac{nk-1}2=\alpha_n$. 

Now let $I$ be any independent set of $\Fn$; we prove that $|I|\le\alpha_n$, with equality only for the sets described in the statement. Two facts will be used repeatedly.
\begin{enumerate}[label=\textup{(\alph*)}]
\item If $1\le\delta\le k-1$ and $I$ contains one of the levels $D_\delta$, $D_{\delta+1}$, then $I$ is disjoint from the other one. Indeed, by Proposition~\ref{prop:levels}(i),(ii), every vertex of each of these levels has a neighbour in the other one.
\item $a_k\le k$. Indeed, if $n=2k$ then $|D_k|=k$, and if $n=2k+1$ then $D_k$ induces a cycle of length $n$ (Proposition~\ref{prop:levels}), whose independence number is $k$.
\end{enumerate}

We distinguish two cases according to the parity of $k$.

\emph{Case 1: $k$ even.} Group the levels into the $k/2$ consecutive pairs $(D_{2j-1},D_{2j})$, $1\le j\le k/2$. For each of them, $a_{2j-1}+a_{2j}\le n$: this follows from Lemma~\ref{lem:pair} if $j<k/2$, or if $j=k/2$ and $n$ is odd (in both cases $|D_{2j}|=n$), and from Lemma~\ref{lem:topeven} if $j=k/2$ and $n$ is even. Since the levels partition $V(\Fn)$, summing over $j$ gives $|I|=\sum_{\delta=1}^k a_\delta\le\frac k2\,n=\alpha_n$. 

Suppose now that $|I|=\alpha_n$. Then all these inequalities are equalities; in particular, $a_{k-1}+a_k=n$. If $n$ is even, Lemma~\ref{lem:topeven} gives $D_{k-1}\subseteq I$ and $I\cap D_k=\emptyset$. If $n$ is odd, Lemma~\ref{lem:pair} shows that $I$ contains one of $D_{k-1}$ and $D_k$ and is disjoint from the other; since $D_k$ is not independent (Proposition~\ref{prop:levels}), we obtain $D_{k-1}\subseteq I$ and $I\cap D_k=\emptyset$. 

\emph{Case 2: $k$ odd.} Group the levels $D_1,\dots,D_{k-1}$ into the $(k-1)/2$ consecutive pairs $(D_{2j-1},D_{2j})$, $1\le j\le (k-1)/2$. For each of them, $a_{2j-1}+a_{2j}\le n$ by Lemma~\ref{lem:pair}, since $2j\le k-1$ and hence $|D_{2j}|=n$. Together with $a_k\le k$, by~(b), summing gives $|I|=\sum_{\delta=1}^k a_\delta\le\frac{k-1}2\,n+k=\alpha_n$.

Suppose now that $|I|=\alpha_n$. Then all these inequalities are equalities; in particular, $W:=I\cap D_k$ has $k$ elements. Hence $W=D_k$ if $n$ is even, since then $|D_k|=k$, and $W$ is a maximum independent set of the cycle induced by $D_k$ if $n$ is odd. If $k=1$, then $I=W$ and we are done. If $k\ge3$, then $W\neq\emptyset$, so $D_{k-1}\not\subseteq I$ by~(a). Since $a_{k-2}+a_{k-1}=n$, Lemma~\ref{lem:pair} shows that $I$ contains one of $D_{k-2}$ and $D_{k-1}$ and is disjoint from the other; therefore $D_{k-2}\subseteq I$ and $I\cap D_{k-1}=\emptyset$.

\emph{Conclusion.} Unless $k=1$, in both cases we have found an odd index $2m+1\in\{k-1,k-2\}$ such that $D_{2m+1}\subseteq I$, $I\cap D_{2m+2}=\emptyset$ and, if $2m+3=k$, $I\cap D_k=W$; moreover, $a_{2j-1}+a_{2j}=n$ for $1\le j\le m$. We show that $D_{2j-1}\subseteq I$ and $I\cap D_{2j}=\emptyset$ for $j=m,m-1,\dots,1$ (there is nothing to prove if $m=0$), by descending induction. Assume that $D_{2j+1}\subseteq I$, which holds for $j=m$. Then $I\cap D_{2j}=\emptyset$ by~(a), and since $a_{2j-1}+a_{2j}=n$, Lemma~\ref{lem:pair} shows that $I$ contains one of $D_{2j-1}$, $D_{2j}$; hence $D_{2j-1}\subseteq I$. Since the levels partition $V(\Fn)$, it follows that $I$ is one of the sets described in the statement.

This proves that the maximum independent sets of $\Fn$ are exactly the sets described in the statement, and hence that $\alpha(\Fn)=\alpha_n$.
\end{proof}

\begin{proof}[Proof of Corollary~\ref{cor:introA}]
For $r\in\Z_n$, the rotation $A\mapsto A+r$ is an automorphism of $\Fn$ that preserves
each level $D_\delta$. Hence each $I_n+r$ is a maximum independent set of $\Fn$, and,
by Theorem~\ref{thm:introA}, all maximum independent sets of $\Fn$ contain the same
levels of odd index below $k$ and can differ only in the block $W\subseteq D_k$.

If $n\not\equiv3\pmod 4$, that is, if $k$ is even or $n$ is even, then $W$ is
determined by Theorem~\ref{thm:introA}, so $\Fn$ has a unique maximum independent set;
by Definition~\ref{def:In} it is $I_n$, and therefore $I_n+r=I_n$ for every $r$.

Now let $n=2k+1$ with $k$ odd, and let $C$ be the cycle induced by $D_k$. By
Theorem~\ref{thm:introA}, the map $W\mapsto D_1\cup D_3\cup\cdots\cup D_{k-2}\cup W$ is a
bijection from the maximum independent sets of $C$, that is, its independent
$k$-subsets, onto the maximum independent sets of $\Fn$. Put $Q_j=\lv{jk}k$ for
$j\in\Z_n$. Since $\gcd(k,n)=1$, the map $j\mapsto Q_j$ is a bijection $\Z_n\to D_k$,
and $Q_j$ is adjacent to $Q_{j\pm1}$ by Proposition~\ref{prop:levels}(iii); hence
$Q_0,Q_1,\dots,Q_{n-1}$ is the cyclic order of $C$.

An independent $k$-subset of $C$ leaves $k$ nonempty gaps between cyclically consecutive
chosen vertices, which contain $n-k=k+1$ vertices in total; so exactly one gap has two
vertices and the others have one. Such a set is determined by the position of its long
gap, and each of the $n$ possible positions gives one; so there are exactly $n$ of them,
and the shift $Q_j\mapsto Q_{j+1}$ permutes them cyclically. The rotation $A\mapsto A+r$
maps $Q_j$ to $\lv{jk+r}k=Q_{j+rk^{-1}}$, where $k^{-1}$ is the inverse of $k$ modulo
$n$; hence the rotations of $\Fn$ induce all the shifts of $C$ and act transitively on
these $n$ sets. Since there are $n$ rotations, each of the $n$ sets is the image of a
given one under exactly one rotation.

Finally, $L_k=\{\lv xk:x\in\cint{k+2}{2k+1}\}$ is one of these sets. It has $k$
elements, and by Proposition~\ref{prop:levels}(iii) two distinct elements $\lv xk$ and
$\lv{x'}k$ of $D_k$ are adjacent if and only if $x'-x=\pm k$, which cannot happen when
$x$ and $x'$ lie in a set of $k$ consecutive residues. Therefore the sets $L_k+r$,
$r\in\Z_n$, are pairwise distinct and are all the maximum independent sets of $C$, and
the sets $I_n+r=D_1\cup\cdots\cup D_{k-2}\cup(L_k+r)$ are exactly the $n$ maximum
independent sets of $\Fn$.
\end{proof}

\begin{remark}\label{rem:paths}
For comparison, we consider the $2$-token graph of a path. Let $P_n$, $n\ge2$, be the path with vertex set
$\{1,\dots,n\}$ and edges $\{i,i+1\}$, $1\le i\le n-1$. For $1\le\delta\le n-1$ let
$D^P_\delta=\bigl\{\{u,v\}:v-u=\delta\bigr\}$, so that $|D^P_\delta|=n-\delta$, and for
an independent set $I$ of $F_2(P_n)$ write $a^P_\delta=|I\cap D^P_\delta|$. The neighbours
of $\{x,x+\delta\}$ in $F_2(P_n)$ are those of the sets $\{x,x+\delta-1\}$,
$\{x+1,x+\delta\}\in D^P_{\delta-1}$ and $\{x-1,x+\delta\}$, $\{x,x+\delta+1\}\in D^P_{\delta+1}$
that are $2$-subsets of $\{1,\dots,n\}$. Hence every edge joins two consecutive levels,
and every level is independent, now with no exception. Moreover, for
$1\le\delta\le n-2$, writing $u_x=\{x,x+\delta\}$ and $w_y=\{y,y+\delta+1\}$, the set
$D^P_\delta\cup D^P_{\delta+1}$ induces the path $u_1w_1u_2w_2\cdots w_{n-\delta-1}u_{n-\delta}$.
This path has odd order, so its unique maximum independent set is $D^P_\delta$. Thus
$a^P_\delta+a^P_{\delta+1}\le n-\delta$ for every independent set $I$, with equality only
if $I\cap(D^P_\delta\cup D^P_{\delta+1})=D^P_\delta$.

Grouping the levels into the pairs $(1,2),(3,4),\dots$, and leaving the level
$D^P_{n-1}=\{\{1,n\}\}$ unpaired when $n$ is even, we get
$|I|\le(n-1)+(n-3)+\cdots=\lfloor n^2/4\rfloor$, where the last term is
$|D^P_{n-1}|=1$ when $n$ is even. This value was already recorded
in~\cite[Section~3]{alba}. Equality forces $I$ to be the union of the levels of odd
index, that is, the set of pairs $\{u,v\}$ with $v-u$ odd. The difference with the cycle
lies in the equality case for consecutive levels. For $C_n$, these levels induce the
cycle $C_{2n}$ of Lemma~\ref{lem:pair}, which has two maximum independent sets, so the
choice must be propagated from the top level. For $P_n$, they induce a path of odd
order, whose maximum independent set is unique, so the choice is forced within each
pair.
\end{remark}

\section{\texorpdfstring{Proof of Theorem~\ref{thm:introB}}{Proof of Theorem 1.3}}\label{sec:main}

In this section we give the proof of Theorem~\ref{thm:introB}. We first explain why the sum map is used instead of a more obvious candidate. A first attempt at a correspondence between maximum independent sets of
$\Kn$ and $\Fn$ is suggested by the isomorphism
$\Fn\cong\bigl((C_n\square C_n)-\Delta\bigr)/S_2$, where $\Delta$ is the
diagonal~\cite{alzaga,aude}. Since $C_n\square C_n$ and $\Kn$ have the same vertex set,
one may send the king at $(x,y)$ to the pair $\pi(x,y)=\{x,y\}$. This works for even
$n$: the set $O\times E$, where $O$ and $E$ are the sets of odd and even elements of
$\Z_n$, respectively, is a maximum independent set of $\Kn$ that $\pi$ maps bijectively onto $I_n$.
However, it fails for small odd values of $n$: a computer search shows that, for
$n\in\{5,7,9\}$, the largest size of a set that avoids the diagonal, is independent in
$\Kn$, and is mapped injectively by $\pi$ to an independent set of $\Fn$ is,
respectively, $4$, $9$ and $16$, whereas $\alpha(\Fn)=5$, $10$ and $18$. The sum map avoids this obstruction. It is obtained from $\pi$ by a change of coordinates: $\psi=\pi\circ\theta$, where $\theta(x,y)=(x,y-x)$. Since $\theta$ is a bijection of $\Z_n\times\Z_n$ that maps $\ell_0$ onto $\Delta$, the map $\psi$ is defined exactly on the squares off $\ell_0$. Since $\theta$ is not an automorphism of $\Kn$, however, it changes which sets of squares are independent, and with the new coordinates a suitable set exists for every $n\ge4$ (Theorem~\ref{thm:main}).

We now fix some definitions and notation. 

Recall that $\ell_0=\{(x,2x):x\in\Z_n\}$ and that
\[
\psi:(\Z_n\times\Z_n)\setminus\ell_0\longrightarrow V(\Fn),\qquad\psi(x,y)=\{x,\,y-x\},
\]
which is well defined since $y-x\neq x$. For every $\{u,v\}\in V(\Fn)$,
\begin{equation}\label{eq:fibre}
\psi^{-1}(\{u,v\})=\{(u,u+v),\ (v,u+v)\},
\end{equation}
so $\psi$ is two-to-one and the two preimages of a pair lie in the same column. Consequently, for $X\subseteq(\Z_n\times\Z_n)\setminus\ell_0$, the restriction $\psi|_X$ is injective if and only if $X$ contains at most one of the two points~\eqref{eq:fibre} of each fibre.

For $c\in\Z_n$, let $\ell_c=\{(x,2x+c):x\in\Z_n\}$ be the \emph{line of slope $2$ and offset $c$}; these lines partition $\Z_n\times\Z_n$, since the offset of $(x,y)$ is $y-2x$. Offsets are elements of $\Z_n$, so conditions on their differences, such as $c'-c\notin\{1,-1\}$ or $c'-c\in\{2,3\}$, are understood in $\Z_n$.

\begin{proposition}\label{prop:lines}
Let $c\in\Z_n\setminus\{0\}$. Then $\psi(x,2x+c)=\{x,x+c\}$ for all $x\in\Z_n$. If moreover $|c|<n/2$, then $\psi|_{\ell_c}\colon\ell_c\to D_{|c|}$ is a bijection, whose inverse is $\{u,v\}\mapsto(u,u+v)$, where the elements of each pair are labelled so that $v=u+c$.
\end{proposition}
\begin{proof}
Since $c\neq0$, we have $2x+c\neq2x$ for every $x\in\Z_n$, so $\ell_c\cap\ell_0=\emptyset$ and $\psi$ is defined on $\ell_c$. For $x\in\Z_n$, $\psi(x,2x+c)=\{x,\,(2x+c)-x\}=\{x,\,x+c\}$, which is the first assertion.

Assume now that $\delta=|c|$ satisfies $0<\delta<n/2$. Since $\dist(x,x+c)=\dist(0,c)=\delta$, the map $\psi$ sends $\ell_c$ into $D_\delta$. Conversely, let $\{u,v\}\in D_\delta$. Then $v-u\in\{c,-c\}$, and $c\neq-c$ because $2c\neq0$ (as $0<|c|<n/2$); hence we may label the elements of the pair, in a unique way, so that $v=u+c$. Define
\[
\chi\colon D_\delta\to\ell_c,\qquad \chi(\{u,v\})=(u,\,u+v)\quad\text{where } v=u+c .
\]
Then $\chi(\{u,v\})=(u,2u+c)\in\ell_c$ and $\psi(\chi(\{u,v\}))=\{u,u+c\}=\{u,v\}$. Moreover, for $x\in\Z_n$, the pair $\psi(x,2x+c)=\{x,x+c\}$ is labelled with $u=x$ and $v=x+c$, so $\chi(\psi(x,2x+c))=(x,2x+c)$. Therefore $\psi|_{\ell_c}$ and $\chi$ are mutually inverse bijections.
\end{proof}

\begin{lemma}\label{lem:attacks}
Let $n\ge4$ and let $p=(x,2x+c)$ and $q=(x',2x'+c')$ be distinct vertices of $\Kn$, and let $\varepsilon=x'-x$. Then $p$ and $q$ attack each other if and only if
\begin{equation}\label{eq:attack}
\varepsilon\in\{0,\pm1\}\quad\text{and}\quad 2\varepsilon+(c'-c)\in\{0,\pm1\}.
\end{equation}
In particular, every line $\ell_c$ is an independent set of $\Kn$.
\end{lemma}

\begin{proof}
Two distinct vertices of $\Kn$ attack each other precisely when their rows and their columns differ by an element of $\{0,\pm1\}$; since the columns differ by $2\varepsilon+(c'-c)$, this is~\eqref{eq:attack}. If $c'=c$, then $\varepsilon\neq0$, so~\eqref{eq:attack} would require $2\varepsilon=\pm2\in\{0,\pm1\}$, which is false for $n\ge4$.
\end{proof}

All three constructions below have the same shape: a family of lines, each one restricted to a prescribed set of rows. The next lemma gives a criterion for such a family to be independent.

\begin{lemma}\label{lem:rows}
Let $n\ge4$ and let $\mathcal E\subseteq\Z_n$ be a set of offsets such that $c'-c\notin\{1,-1\}$ for all distinct $c,c'\in\mathcal E$. For each $c\in\mathcal E$ choose $\mathcal B_c\subseteq\Z_n$ and let
\[
S=\bigcup_{c\in\mathcal E}\{(x,\,2x+c):x\in \mathcal B_c\}.
\]
If $\mathcal B_c\cap(\mathcal B_{c'}+1)=\emptyset$ whenever $c,c'\in\mathcal E$ satisfy $c'-c\in\{2,3\}$, then $S$ is an independent set of $\Kn$.
\end{lemma}

\begin{proof}
Let $p=(x,2x+c)$ and $q=(x',2x'+c')$ be distinct vertices of $S$, with $c,c'\in\mathcal E$, $x\in \mathcal B_c$ and $x'\in \mathcal B_{c'}$, and put $\varepsilon=x'-x$. By~\eqref{eq:attack}, $p$ and $q$ attack each other only if $\varepsilon\in\{0,\pm1\}$ and $2\varepsilon+(c'-c)\in\{0,\pm1\}$. If $\varepsilon=0$, then $c'-c\in\{0,\pm1\}$; the values $\pm1$ are excluded by hypothesis when $c\neq c'$, so $c'=c$ and $p=q$, a contradiction. If $\varepsilon=1$, then $c-c'\in\{1,2,3\}$, hence $c-c'\in\{2,3\}$; but $x\in \mathcal B_c$ and $x+1=x'\in \mathcal B_{c'}$, contradicting $\mathcal B_{c'}\cap(\mathcal B_c+1)=\emptyset$. If $\varepsilon=-1$, then $c'-c\in\{1,2,3\}$, hence $c'-c\in\{2,3\}$; but $x'=x-1\in \mathcal B_{c'}$ gives $x\in \mathcal B_{c'}+1$, and $x\in \mathcal B_c$, contradicting $\mathcal B_c\cap(\mathcal B_{c'}+1)=\emptyset$.
\end{proof}

We shall also need to know when two points lying on lines have the same image under $\psi$.

\begin{lemma}\label{lem:collide}
Let $c,c'\in\Z_n\setminus\{0\}$ and $x,x'\in\Z_n$ with $(x,c)\neq(x',c')$. Then $\psi(x,2x+c)=\psi(x',2x'+c')$ if and only if $c'=-c$ and $x'=x+c$.
\end{lemma}

\begin{proof}
By Proposition~\ref{prop:lines}, the condition reads $\{x,x+c\}=\{x',x'+c'\}$. If $x'=x$, it gives $c'=c$, which is excluded; hence $x'=x+c$ and $x'+c'=x$, that is, $c'=-c$. The converse is clear.
\end{proof}

\subsection{The construction}
Proposition~\ref{prop:lines} and Lemma~\ref{lem:attacks} suggest how to build $\Sn$:
for $1\le|c|<n/2$, the line $\ell_c$ is an independent set of $\Kn$ that $\psi$ maps
bijectively onto the level $D_{|c|}$. More generally, one may restrict each line to a
set of rows, and Lemma~\ref{lem:rows} gives a sufficient condition for the resulting set
to be independent. Three cases occur, according to whether $n$ is even,
$n\equiv1\pmod4$ or $n\equiv3\pmod4$.

If $n$ is even, then the elements of $\Z_n$ have a well-defined parity, and it is enough
to take the same set of rows, namely the odd ones, for \emph{every} odd offset (Case~1 of
the proof of Theorem~\ref{thm:main}). If $n\equiv1\pmod4$, then $I_n$ is the union of the complete levels $D_\delta$ with $\delta$ odd and $\delta\le k-1$, so it suffices to choose, for each such $\delta$, the offset $c\in\{\pm\delta\}$ with $c\equiv1\pmod4$ and to keep the line $\ell_c$ whole. The chosen offsets differ by nonzero multiples of~$4$ of absolute value at most $n-5$, so no two of them differ by $\pm1$, $2$ or $3$ modulo $n$; hence the hypotheses of Lemma~\ref{lem:rows} hold, the second one vacuously (Case~2 of the proof). If $n\equiv3\pmod4$, then $I_n$ contains only $k$ of the $n$ elements of the top level $D_k$, because $D_k$ is not independent; these $k$ elements are obtained from $k$ squares of the line $\ell_{-k}$. For each odd $\delta<k$, the whole level $D_\delta$ is needed, but now it is obtained from parts of both lines $\ell_\delta$ and $\ell_{-\delta}$, namely from the sets $R_t$ with $c_t=\pm\delta$ defined below: $\psi$ maps the squares chosen on one line onto some pairs of $D_\delta$, and those chosen on the other onto the remaining pairs.

We now introduce the sets of offsets and rows used in the three cases; the construction itself is carried out in the proof of Theorem~\ref{thm:main}.

For the first two cases, let
\begin{align*}
 O&=\{1,3,\dots,n-1\}\subseteq\Z_n \hspace{3.4cm} (n \text{ even}),\\
\mathcal C_k&=\{c\in\Z:\ |c|\le k-1,\ c\equiv1\!\!\pmod 4\} \hspace{1.0cm}  (k \text{ even}).
\end{align*}

For $n=2k+1$ with $k$ odd and $0\le t\le k-1$, let 
\[
c_t=2t-k
\] 
and
\[
J_t=\cint{k-t+1}{2k-t},\qquad
B_t=\begin{cases}
\cint 1k, & t=0,\\
J_t, & t\text{ odd},\\
\Z_n\setminus J_t, & t\ge2\text{ even}.
\end{cases}
\]
The set $B_t$ is the set of rows used on the line $\ell_{c_t}$: the points $(x,2x+c_t)$ with $x\in B_t$ form the subset $R_t$ of $\ell_{c_t}$ that appears in Theorem~\ref{thm:main}. Note that $|J_t|=k$, that $\Z_n\setminus J_t=\cint{2k-t+1}{3k-t+1}$ has $k+1$ elements, and that $B_0\subseteq\Z_n\setminus J_0=\{n\}\cup\cint 1k$. The role of the sets $B_t$ is to satisfy the two conditions in the following lemma, the first for the independence of $\Sn$ and the second for the injectivity of $\psi$ on $\Sn$. These two conditions are the reason why $B_t$ alternates between $J_t$ and its complement. Here and below, the indices $t$ of $B_t$ are read modulo $k$.

\begin{lemma}\label{lem:chain}
Let $n=2k+1$ with $k\ge3$ odd. Then:
\begin{enumerate}[label=\textup{(\roman*)}]
\item $B_t\cap(B_{t+1}+1)=\emptyset$ for every $t\in\Z_k$;
\item $B_t\sqcup\bigl(B_{k-t}+(k-2t)\bigr)=\Z_n$ for $1\le t<k/2$.
\end{enumerate}
\end{lemma}

\begin{proof}
We first note that $J_{t+1}+1=J_t$ for $0\le t\le k-2$. Indeed, $J_{t+1}+1=\cint{k-t}{2k-t-1}+1=\cint{k-t+1}{2k-t}=J_t$. Therefore, by iteration we have 
\begin{equation}\label{eq:Jshift}
J_{t+m}+m=J_t\quad(0\le t\le t+m\le k-1).
\end{equation}
(i) It is enough to show that $B_{t+1}+1\subseteq\Z_n\setminus B_t$. If $t$ is odd (so $1\le t\le k-2$ and $t+1\ge2$ is even), then $B_{t+1}+1=(\Z_n\setminus J_{t+1})+1=\Z_n\setminus(J_{t+1}+1)=\Z_n\setminus J_t=\Z_n\setminus B_t$. If $t$ is even with $2\le t\le k-3$, then $B_{t+1}+1=J_{t+1}+1=J_t=\Z_n\setminus B_t$. If $t=0$, then $B_1+1=J_1+1=J_0$, which is disjoint from $B_0\subseteq\Z_n\setminus J_0$. Finally, if $t=k-1$, which is even, then $t+1\equiv0\pmod k$, so $B_{t+1}=B_0$ and $B_0+1=\cint 2{k+1}=J_{k-1}=\Z_n\setminus B_{k-1}$.

(ii) Put $s=k-t$ and $\delta=k-2t$, so that $1\le t<s\le k-1$, $s-\delta=t$ and $t+s=k$ is odd, whence $t$ and $s$ have different parities. By~\eqref{eq:Jshift}, $J_s+\delta=J_{s-\delta}=J_t$. Hence, if $t$ is odd then $B_t=J_t$ and $B_s+\delta=(\Z_n\setminus J_s)+\delta=\Z_n\setminus(J_s+\delta)=\Z_n\setminus J_t$, while if $t$ is even then $B_t=\Z_n\setminus J_t$ and $B_s+\delta=J_s+\delta=J_t$. In both cases $B_t$ and $B_s+\delta$ partition $\Z_n$.
\end{proof}

\begin{theorem}\label{thm:main}
Let $n\ge4$, $k=\lfloor n/2\rfloor$, and let
\[
\Sn=\begin{cases}
O\times O, & n\text{ even},\\[4pt]
\displaystyle\bigcup_{c\in \mathcal C_k}\ell_c, & n\equiv1\!\!\pmod 4,\\[10pt]
\displaystyle\bigcup_{t=0}^{k-1}R_t,\quad\text{where } R_t=\{(x,2x+c_t):x\in B_t\}, & n\equiv3\!\!\pmod 4.
\end{cases}
\]
Then $\Sn$ is a maximum independent set of $\Kn$, $\Sn\cap\ell_0=\emptyset$, and $\psi$ maps $\Sn$ bijectively onto $I_n$. In particular, Theorem~\ref{thm:introB} holds.
\end{theorem}

\begin{proof}
In each of the three cases, we check that $\Sn\cap\ell_0=\emptyset$, that $\Sn$ is independent, and that $\psi$ maps $\Sn$ bijectively onto $I_n$. This suffices: $\Sn$ is then a maximum independent set, because $|\Sn|=|I_n|=\alpha_n=\alpha(\Kn)$ by Theorem~\ref{thm:introA} and~\eqref{eq:alphas}. Independence always follows from Lemma~\ref{lem:rows}. To prove that $\psi|_{\Sn}$ is a bijection onto $I_n$, in Case~2 we use directly that $\psi$ maps each line $\ell_c$ bijectively onto $D_{|c|}$ (Proposition~\ref{prop:lines}). In Cases~1 and~3 we show instead that $\psi(\Sn)\subseteq I_n$, that $\psi|_{\Sn}$ is injective and that $|\Sn|=\alpha_n$. Since $|I_n|=\alpha_n$ (this is computed in the proof of Theorem~\ref{thm:introA}; see also~\cite{alba}), an injection from $\Sn$ into $I_n$ is then a bijection. 

\medskip\noindent\textbf{Case 1: $n$ even.}
By Remark~\ref{rem:In}, the set $I_n$ consists of the pairs with one odd and one even element.

\emph{Bijectivity.} Let $(x,y)\in O\times O$. Then $2x$ is even and $y$ is odd, so $(x,y)\notin\ell_0$; and $\psi(x,y)=\{x,y-x\}$ has one odd and one even element, so it lies in $I_n$. The other preimage $(y-x,y)$ of this pair has even row, so it does not lie in $O\times O$; hence $\psi|_{O\times O}$ is injective. Finally $|O\times O|=k^2=\alpha_n$.

\emph{Independence.} A pair $(x,y)$ with $x$ odd lies in $O\times O$ if and only if the offset $c=y-2x$ is odd, so $O\times O$ is the union of the sets $\{(x,2x+c):x\in O\}$ over the odd $c\in\Z_n$; that is, $\mathcal B_c=O$ for every odd offset. Two odd offsets differ by an even amount, which is never $\pm1$ because $n$ is even, and never $3$; and if $c'-c=2$, then $\mathcal B_c\cap(\mathcal B_{c'}+1)=O\cap(O+1)=\emptyset$, since $O+1$ is the set of even elements of $\Z_n$. Lemma~\ref{lem:rows} now gives that $\Sn$ is independent.

\medskip\noindent\textbf{Case 2: $n\equiv1\pmod 4$.}
Here $n$ is odd and $k$ is even. The odd integers in $[-(k-1),k-1]$ are
$\pm1,\pm3,\dots,\pm(k-1)$. For each odd $\delta$, the integers $\delta$ and $-\delta$
are odd, hence each is congruent to $1$ or $3$ modulo $4$, and they are not congruent to
each other because $\delta-(-\delta)=2\delta\equiv2\pmod4$; so exactly one of them is
congruent to $1$ modulo $4$. Hence $c\mapsto|c|$ is a bijection from $\mathcal C_k$ onto
$\{1,3,\dots,k-1\}$, and $|\mathcal C_k|=k/2$. Moreover, the elements of $\mathcal C_k$
are nonzero modulo $n$, being odd with $|c|<n$, and pairwise incongruent modulo $n$,
since $|c-c'|\le2k-2<n$; thus we may regard $\mathcal C_k$ as a subset of
$\Z_n\setminus\{0\}$ with $k/2$ elements.

\emph{Bijectivity.} Since the offsets in $\mathcal C_k$ are nonzero, $\Sn\cap\ell_0=\emptyset$. By Proposition~\ref{prop:lines}, $\psi$ maps each line $\ell_c$, $c\in\mathcal C_k$, bijectively onto $D_{|c|}$. As the lines $\ell_c$ are pairwise disjoint and $c\mapsto|c|$ is a bijection from $\mathcal C_k$ onto $\{1,3,\dots,k-1\}$, the map $\psi$ sends $\Sn$ bijectively onto $D_1\cup D_3\cup\cdots\cup D_{k-1}=I_n$.

\emph{Independence.} Here $\mathcal B_c=\Z_n$ for $c\in\mathcal C_k$. If $c$ and $c'$ are
distinct elements of $\mathcal C_k$, then $c'-c=4j$ for some integer $j\neq0$, because
both are congruent to $1$ modulo $4$. Moreover, $|c'-c|\le2k-2$, and since $k$ is even,
$2k-2\equiv2\pmod4$; hence $|c'-c|\le2k-4$, that is, $1\le|j|\le k/2-1$. Therefore
$4\le|c'-c|\le2k-4=n-5$, and $c'-c$ is congruent modulo $n$ to no element of
$\{0,\pm1,2,3\}$. Thus the hypotheses of Lemma~\ref{lem:rows} hold, the second one
vacuously, and $\Sn$ is independent.

\medskip\noindent\textbf{Case 3: $n\equiv3\pmod 4$.}
Here $n$ is odd and $k$ is odd, so $k\ge3$ and $n\ge7$. The offsets $c_0=-k,\ c_1=-k+2,\ \dots,\ c_{k-1}=k-2$ are odd, nonzero and pairwise incongruent modulo $n$, because $|c_s-c_t|\le2k-2<n$. In particular $\Sn\cap\ell_0=\emptyset$. Moreover, $R_t\subseteq\ell_{c_t}$ for each $t$, so the sets $R_t$ lie on pairwise distinct lines; hence they are pairwise disjoint, and $\Sn\cap\ell_{c_t}=R_t$ for every $t$. 

\emph{Image.} For $x\in B_0=\cint1k$ we have $\psi(x,2x-k)=\{x-k,x\}=\lv{x-k}k$ with $x-k\in\cint{k+2}{2k+1}$; hence $\psi(R_0)\subseteq L_k$. For $1\le t\le k-1$, the integer $|c_t|=|2t-k|$ is odd and at most $k-2$, so $\psi(R_t)\subseteq D_{|c_t|}\subseteq I_n$ by Proposition~\ref{prop:lines}. Therefore $\psi(\Sn)\subseteq I_n$.

\emph{Injectivity.} Let $p\neq p'$ be points of $\Sn$ with $\psi(p)=\psi(p')$, say $p\in R_t$ and $p'\in R_s$. Since $R_t\subseteq\ell_{c_t}$ and $R_s\subseteq\ell_{c_s}$, Lemma~\ref{lem:collide} gives $c_s\equiv-c_t\pmod n$, and, writing $p=(x,2x+c_t)$, also $p'=(x+c_t,\,2(x+c_t)+c_s)$. Now $c_s+c_t=2(s+t-k)$ and $|c_s+c_t|\le2k<n$, so $s+t=k$. As $s,t\le k-1$, neither of them is $0$, and since $k$ is odd, $s\neq t$. Exchanging $p$ and $p'$ if necessary, we may assume $1\le t<k/2$ and $s=k-t$. Let $\delta=k-2t>0$, so that $c_t=-\delta$ and $c_s=\delta$. Then $p=(x,2x-\delta)$ and $p'=(x-\delta,\,2(x-\delta)+\delta)$, with $x\in B_t$ and $x-\delta\in B_s$; that is, $x\in B_t\cap(B_s+\delta)$. This set is empty by Lemma~\ref{lem:chain}(ii), a contradiction. Hence $\psi|_{\Sn}$ is injective.

\emph{Cardinality.} Among $1,\dots,k-1$ there are $(k-1)/2$ odd and $(k-1)/2$ even values of $t$, with $|B_t|=k$ and $|B_t|=k+1$ respectively, while $|B_0|=k$. Hence
\[
|\Sn|=\sum_{t=0}^{k-1}|B_t|=k+\tfrac{k-1}2\,k+\tfrac{k-1}2\,(k+1)=k+\tfrac{(k-1)n}2=\tfrac{nk-1}2=\alpha_n .
\]

\emph{Independence.} Here $\mathcal B_{c_t}=B_t$. The differences $c_s-c_t=2(s-t)$ are even and satisfy $|c_s-c_t|\le2k-2=n-3$, so none of them is congruent to $\pm1$ modulo $n$. Moreover $c_s-c_t\equiv2$ forces $s=t+1$, and $c_s-c_t\equiv3$ forces $2(s-t)=3-n=-(2k-2)$, that is, $t=k-1$ and $s=0$. In both cases the condition required by Lemma~\ref{lem:rows} reads $B_t\cap(B_{t+1}+1)=\emptyset$, with the indices of the $B_t$ read modulo $k$, and this is Lemma~\ref{lem:chain}(i). Hence $\Sn$ is independent. Example~\ref{ex:n7} works out the case $n=7$ in full.
\end{proof}

\subsection{\texorpdfstring{The map $\Phi_n$}{The map Phin}}
\begin{definition}\label{def:Phi}
For $n\ge4$, we denote by $\Phi_n\colon I_n\to\Sn$ the inverse of the bijection $\psi|_{\Sn}$ of Theorem~\ref{thm:main}.
\end{definition}

By~\eqref{eq:fibre}, $\Phi_n(\{u,v\})=(r,u+v)$ for a token $r\in\{u,v\}$, and the next result gives the rule that selects it. Whenever $\delta=\dist(u,v)<n/2$, the equality $\{u,v\}=\{u,u+\delta\}$ determines $u$ and $v$ uniquely, and we label them so that $v=u+\delta$. For odd $n$ every distance occurring in $I_n$ is at most $k<n/2$; for even $n$ no labelling is needed.

\begin{proposition}\label{prop:Phi}
Let $\{u,v\}\in I_n$ and $\delta=\dist(u,v)$. Then $\Phi_n(\{u,v\})=(r,u+v)$, where:
\begin{enumerate}[label=\textup{(\roman*)}]
\item if $n$ is even, then $r$ is the odd element of $\{u,v\}$;
\item if $n\equiv1\pmod4$, then $r=u$ when $\delta\equiv1\pmod4$ and $r=v$ when $\delta\equiv3\pmod4$;
\item if $n\equiv3\pmod4$, then $r=v$ when $v\in B_{(k-\delta)/2}$, and $r=u$ otherwise.
\end{enumerate}
\end{proposition}

\begin{proof}
Exactly one of the two points $(u,u+v)$, $(v,u+v)$ lies in $\Sn$, so it suffices to check that $(r,u+v)\in\Sn$. (i) Exactly one of $u,v$ is odd (Remark~\ref{rem:In}), so both $r$ and $u+v$ are odd, and $(r,u+v)\in O\times O$. (ii) If $\delta\equiv1\pmod 4$, then $(u,u+v)=(u,2u+\delta)\in\ell_\delta$ and $\delta\in\mathcal C_k$; if $\delta\equiv3\pmod4$, then $(v,u+v)=(v,2v-\delta)\in\ell_{-\delta}$ and $-\delta\in\mathcal C_k$. (iii) Put $t=(k-\delta)/2$, so that $c_t=-\delta$ and $(v,u+v)=(v,2v+c_t)$. Since $\Sn\cap\ell_{c_t}=R_t$, this point lies in $\Sn$ if and only if $v\in B_t$; when $v\notin B_t$, the point of the fibre lying in $\Sn$ is therefore $(u,u+v)$.
\end{proof}

\begin{figure}[htbp]
\centering
\definecolor{dist1}{HTML}{A9BAD6}
\definecolor{dist3}{HTML}{E9C98D}
\begin{tikzpicture}[scale=0.5,baseline=(current bounding box.north)]
\fill[dist1] (2,6) rectangle (3,7);
\node[font=\tiny] at (2.5,6.5) {$1,2$};
\fill[dist3] (5,6) rectangle (6,7);
\node[font=\tiny] at (5.5,6.5) {$1,5$};
\fill[dist3] (0,5) rectangle (1,6);
\node[font=\tiny] at (0.5,5.5) {$2,6$};
\fill[dist3] (2,4) rectangle (3,5);
\node[font=\tiny] at (2.5,4.5) {$3,7$};
\fill[dist1] (4,4) rectangle (5,5);
\node[font=\tiny] at (4.5,4.5) {$2,3$};
\fill[dist1] (6,3) rectangle (7,4);
\node[font=\tiny] at (6.5,3.5) {$3,4$};
\fill[dist1] (1,2) rectangle (2,3);
\node[font=\tiny] at (1.5,2.5) {$4,5$};
\fill[dist1] (3,2) rectangle (4,3);
\node[font=\tiny] at (3.5,2.5) {$5,6$};
\fill[dist1] (5,1) rectangle (6,2);
\node[font=\tiny] at (5.5,1.5) {$6,7$};
\fill[dist1] (0,0) rectangle (1,1);
\node[font=\tiny] at (0.5,0.5) {$1,7$};
\draw[step=1,gray!70,thin] (0,0) grid (7,7);
\draw[thick] (0,0) rectangle (7,7);
\node[font=\scriptsize] at (0.5,7.5) {$1$};
\node[font=\scriptsize] at (-0.5,6.5) {$1$};
\node[font=\scriptsize] at (1.5,7.5) {$2$};
\node[font=\scriptsize] at (-0.5,5.5) {$2$};
\node[font=\scriptsize] at (2.5,7.5) {$3$};
\node[font=\scriptsize] at (-0.5,4.5) {$3$};
\node[font=\scriptsize] at (3.5,7.5) {$4$};
\node[font=\scriptsize] at (-0.5,3.5) {$4$};
\node[font=\scriptsize] at (4.5,7.5) {$5$};
\node[font=\scriptsize] at (-0.5,2.5) {$5$};
\node[font=\scriptsize] at (5.5,7.5) {$6$};
\node[font=\scriptsize] at (-0.5,1.5) {$6$};
\node[font=\scriptsize] at (6.5,7.5) {$7$};
\node[font=\scriptsize] at (-0.5,0.5) {$7$};
\node[font=\small] at (3.5,-0.7) {$n=7$};
\end{tikzpicture}
\hspace{0.9cm}
\begin{tikzpicture}[scale=0.5,baseline=(current bounding box.north)]
\fill[dist1] (2,8) rectangle (3,9);
\node[font=\tiny] at (2.5,8.5) {$1,2$};
\fill[dist3] (7,8) rectangle (8,9);
\node[font=\tiny] at (7.5,8.5) {$1,7$};
\fill[dist3] (0,7) rectangle (1,8);
\node[font=\tiny] at (0.5,7.5) {$2,8$};
\fill[dist1] (4,7) rectangle (5,8);
\node[font=\tiny] at (4.5,7.5) {$2,3$};
\fill[dist3] (2,6) rectangle (3,7);
\node[font=\tiny] at (2.5,6.5) {$3,9$};
\fill[dist1] (6,6) rectangle (7,7);
\node[font=\tiny] at (6.5,6.5) {$3,4$};
\fill[dist3] (4,5) rectangle (5,6);
\node[font=\tiny] at (4.5,5.5) {$1,4$};
\fill[dist1] (8,5) rectangle (9,6);
\node[font=\tiny] at (8.5,5.5) {$4,5$};
\fill[dist1] (1,4) rectangle (2,5);
\node[font=\tiny] at (1.5,4.5) {$5,6$};
\fill[dist3] (6,4) rectangle (7,5);
\node[font=\tiny] at (6.5,4.5) {$2,5$};
\fill[dist1] (3,3) rectangle (4,4);
\node[font=\tiny] at (3.5,3.5) {$6,7$};
\fill[dist3] (8,3) rectangle (9,4);
\node[font=\tiny] at (8.5,3.5) {$3,6$};
\fill[dist3] (1,2) rectangle (2,3);
\node[font=\tiny] at (1.5,2.5) {$4,7$};
\fill[dist1] (5,2) rectangle (6,3);
\node[font=\tiny] at (5.5,2.5) {$7,8$};
\fill[dist3] (3,1) rectangle (4,2);
\node[font=\tiny] at (3.5,1.5) {$5,8$};
\fill[dist1] (7,1) rectangle (8,2);
\node[font=\tiny] at (7.5,1.5) {$8,9$};
\fill[dist1] (0,0) rectangle (1,1);
\node[font=\tiny] at (0.5,0.5) {$1,9$};
\fill[dist3] (5,0) rectangle (6,1);
\node[font=\tiny] at (5.5,0.5) {$6,9$};
\draw[step=1,gray!70,thin] (0,0) grid (9,9);
\draw[thick] (0,0) rectangle (9,9);
\node[font=\scriptsize] at (0.5,9.4) {$1$};
\node[font=\scriptsize] at (-0.5,8.5) {$1$};
\node[font=\scriptsize] at (1.5,9.4) {$2$};
\node[font=\scriptsize] at (-0.5,7.5) {$2$};
\node[font=\scriptsize] at (2.5,9.4) {$3$};
\node[font=\scriptsize] at (-0.5,6.5) {$3$};
\node[font=\scriptsize] at (3.5,9.4) {$4$};
\node[font=\scriptsize] at (-0.5,5.5) {$4$};
\node[font=\scriptsize] at (4.5,9.4) {$5$};
\node[font=\scriptsize] at (-0.5,4.5) {$5$};
\node[font=\scriptsize] at (5.5,9.4) {$6$};
\node[font=\scriptsize] at (-0.5,3.5) {$6$};
\node[font=\scriptsize] at (6.5,9.4) {$7$};
\node[font=\scriptsize] at (-0.5,2.5) {$7$};
\node[font=\scriptsize] at (7.5,9.4) {$8$};
\node[font=\scriptsize] at (-0.5,1.5) {$8$};
\node[font=\scriptsize] at (8.5,9.4) {$9$};
\node[font=\scriptsize] at (-0.5,0.5) {$9$};
\node[font=\small] at (4.5,-0.7) {$n=9$};
\end{tikzpicture}
\caption{The sets $\Sset{7}$ (left) and $\Sset{9}$ (right); rows $x$ increase downwards and columns $y$ increase to the right. Each square $(x,y)\in \Sn$ is labelled with $\psi(x,y)=\{x,y-x\}\in I_n$; blue and orange squares correspond to distances $1$ and $3$, respectively.}
\label{fig:examples}
\end{figure}
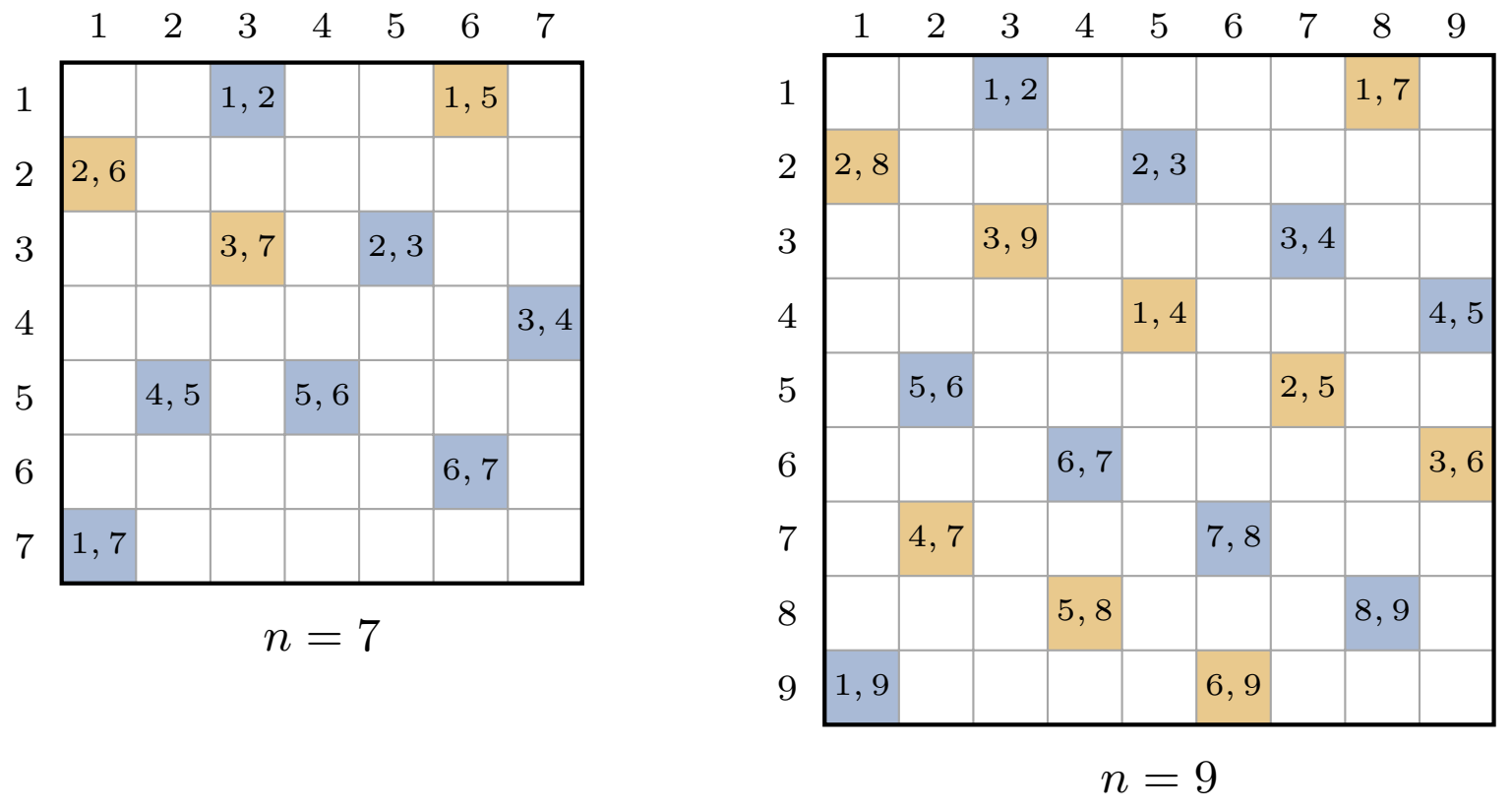

\begin{example}\label{ex:n7}
For $n=7$ we have $k=3$, $c_0=-3$, $c_1=-1$, $c_2=1$, $B_0=\{1,2,3\}$, $B_1=\{3,4,5\}$ and $B_2=\{5,6,7,1\}$, so that
\[
\Sset{7}=\{(1,6),(2,1),(3,3)\}\cup\{(3,5),(4,7),(5,2)\}\cup\{(5,4),(6,6),(7,1),(1,3)\}
\]
and $\Phi_7$ is given by
\[
\begin{array}{c@{\ \mapsto\ }c@{\qquad}c@{\ \mapsto\ }c@{\qquad}c@{\ \mapsto\ }c@{\qquad}c@{\ \mapsto\ }c@{\qquad}c@{\ \mapsto\ }c}
\{1,5\}&(1,6) & \{2,6\}&(2,1) & \{3,7\}&(3,3) & \{2,3\}&(3,5) & \{3,4\}&(4,7)\\[2pt]
\{4,5\}&(5,2) & \{5,6\}&(5,4) & \{6,7\}&(6,6) & \{7,1\}&(7,1) & \{1,2\}&(1,3).
\end{array}
\]
For $n=9$ we have $k=4$ and $\mathcal C_4=\{-3,1\}$, so $\Sset9=\ell_1\cup\ell_{-3}$; this set is displayed in Figure~\ref{fig:examples}. For $n=8$, by contrast, $\Sset8=O\times O$ with $O=\{1,3,5,7\}$.
\end{example}

\subsection{All the maximum independent sets}
By combining Theorem~\ref{thm:main} with the classification in Section~\ref{sec:classification}, we conclude that every maximum independent set of $\Fn$ arises, via the sum map, from a maximum independent set of $\Kn$.

\begin{corollary}\label{cor:all}
Let $n\ge4$. For every maximum independent set $I$ of $\Fn$ there is a maximum independent set $S$ of $\Kn$ with $S\cap\ell_0=\emptyset$ such that $\psi$ restricts to a bijection from $S$ onto $I$.
\end{corollary}

\begin{proof}
For $r\in\Z_n$ the translation $\tau_r(x,y)=(x+r,y+2r)$ is an automorphism of $\Kn$ and
\[
\psi(\tau_r(x,y))=\{x+r,\ y-x+r\}=\psi(x,y)+r .
\]
Moreover $\tau_r$ preserves each line $\ell_c$, since the offset of $\tau_r(x,y)$ is $(y+2r)-2(x+r)=y-2x$; in particular $\tau_r(\Sn)\cap\ell_0=\emptyset$, so $\psi$ is defined on $\tau_r(\Sn)$ and maps it bijectively onto $I_n+r$, and $\tau_r(\Sn)$ is again a maximum independent set of $\Kn$. By Corollary~\ref{cor:introA}, every maximum independent set of $\Fn$ is of the form $I=I_n+r$ for some $r\in\Z_n$, and then $S=\tau_r(\Sn)$ has the required properties. 
\end{proof}

\section{More than two tokens}\label{sec:three}
It is natural to ask whether the coincidence~\eqref{eq:alphas} persists for $m\ge3$ tokens,
with $\Kn$ replaced by the graph $C_n^{\boxtimes m}=C_n\boxtimes\cdots\boxtimes C_n$ ($m$ factors) of non-attacking kings on the
$m$-dimensional torus of side $n$. It does not, and the reason can be seen by comparing densities. Below, $\nu(G)$ denotes the matching number of $G$.

\begin{proposition}\label{prop:three}
For every $n\ge4$ and every $3\le m\le n-1$ we have $\alpha(F_m(C_n))<\alpha(C_n^{\boxtimes m})$.
\end{proposition}

Since $A\mapsto\Z_n\setminus A$ is an isomorphism $F_m(C_n)\cong F_{n-m}(C_n)$, for $m$ close to $n$ the inequality is immediate (for instance, $F_{n-1}(C_n)\cong C_n$); the proof below, however, treats all values of $m$ uniformly.

\begin{proof}
We compare a lower bound for $\alpha(C_n^{\boxtimes m})$ with an upper bound for
$\alpha(F_m(C_n))$. Recall that $k=\lfloor n/2\rfloor\ge2$ and $\alpha_n=\alpha(\Kn)$.

\emph{Lower bound.} If $I$ and $J$ are independent sets of graphs $G$ and $H$, then
$I\times J$ is an independent set of $G\boxtimes H$. Indeed, two distinct vertices
$(g,h)$ and $(g',h')$ of $I\times J$ can be adjacent only if $g=g'$ or
$gg'\in E(G)$, and the independence of $I$ forces $g=g'$; likewise $h=h'$, a
contradiction. Applying this repeatedly, starting from $\Kn$ and adding one factor
$C_n$ at a time, and using~\eqref{eq:alphas} and $\alpha(C_n)=k$, we obtain
\[
\alpha(C_n^{\boxtimes m})\ \ge\ \alpha_n\,k^{m-2}=:L_{n,m}.
\]

\emph{Upper bound.} Since $\nu(C_n)=k$, \cite[Thm.~1.1]{alba}, applied with $m$
tokens, gives
\[
\nu(F_m(C_n))=\tfrac12\tbinom nm\quad\text{if $n$ is even and $m$ is odd},\qquad
\nu(F_m(C_n))\ge\tfrac12\Bigl(\tbinom nm-\tbinom k{\lfloor m/2\rfloor}\Bigr)\quad\text{otherwise}.
\]
An independent set contains at most one end of each edge of a matching, so
$\alpha(F_m(C_n))\le\binom nm-\nu(F_m(C_n))$. Hence
\[
\alpha(F_m(C_n))\ \le\ \tfrac12\Bigl(\tbinom nm+\varepsilon\tbinom k{\lfloor m/2\rfloor}\Bigr)=:U_{n,m},
\]
where $\varepsilon=0$ if $n$ is even and $m$ is odd, and $\varepsilon=1$ otherwise.
It remains to show that $U_{n,m}<L_{n,m}$.

\emph{The case $m=3$.} Here $U_{n,3}=\tfrac12\bigl(\binom n3+\varepsilon k\bigr)$ and
$L_{n,3}=\alpha_nk$.
\begin{itemize}
\item If $n=2k$, then $\varepsilon=0$, $U_{n,3}=\tfrac13k(k-1)(2k-1)$ and $L_{n,3}=k^3$, so
$L_{n,3}-U_{n,3}=\tfrac13k(k^2+3k-1)>0$.
\item If $n=2k+1$, then $\varepsilon=1$, $U_{n,3}=\tfrac13k(2k^2+1)$ and, since
$\alpha_n\ge\frac12(nk-1)$, $L_{n,3}\ge\tfrac12k(2k^2+k-1)$; so
$L_{n,3}-U_{n,3}\ge\tfrac16k(k-1)(2k+5)>0$.
\end{itemize}

\emph{The case $m\ge4$.} We use the estimates
\[
\tbinom nm\le\frac{n(n-1)^{m-1}}{m!},\qquad \frac{2^{m-1}}{m!}\le\frac13,\qquad
\tbinom k{\lfloor m/2\rfloor}\le k^{\lfloor m/2\rfloor}\le k^{m-2}.
\]
\begin{itemize}
\item If $n=2k$, then $\binom nm<(2k)^m/m!\le\tfrac23k^m$ and $L_{n,m}=k^m$, so
$2U_{n,m}<\tfrac23k^m+k^{m-2}<2k^m=2L_{n,m}$.
\item If $n=2k+1$, then $\binom nm\le(2k+1)(2k)^{m-1}/m!\le\tfrac13(2k+1)k^{m-1}$ and
$2\alpha_n\ge nk-1$, so
\[
2L_{n,m}-2U_{n,m}\ \ge\ \bigl((2k+1)k-1\bigr)k^{m-2}-\tfrac13(2k+1)k^{m-1}-k^{m-2}
=\bigl(\tfrac23(2k+1)k-2\bigr)k^{m-2}>0.
\]
\end{itemize}

In all cases, $\alpha(F_m(C_n))\le U_{n,m}<L_{n,m}\le\alpha(C_n^{\boxtimes m})$.
\end{proof}

\begin{remark}\label{rem:three}
The proof of Proposition~\ref{prop:three} also suggests why two tokens are special.
For fixed $m$ and large $n$, the bounds $U_{n,m}$ and $L_{n,m}$ behave like
$n^m/(2\cdot m!)$ and $n^m/2^m$, respectively, and the leading constants agree
precisely when $2\cdot m!=2^m$, that is, when $m\le2$.
\end{remark}

\section*{Acknowledgements}
The author thanks Paloma Jim\'enez-Sep\'ulveda, then an undergraduate student and coauthor of~\cite{JR}, for useful suggestions on the search for a bijection. In the final stage of the project, the author used Anthropic's Claude to refine the wording and grammar, assist with the computational verification of the bijection via concrete examples for small values of $n$, and help check and refine certain proof arguments. All AI-generated code, logical steps, and mathematical arguments were rigorously audited and verified by the author, who assumes full responsibility for the integrity and accuracy of the final work.

\end{document}